\documentclass{amsart}

\usepackage{amssymb,latexsym,amsfonts,amsmath}
\usepackage{graphicx,psfrag,epsfig,epsf}
\usepackage{eso-pic}
\usepackage{graphicx}
\usepackage{color}

\usepackage{tikz}

\usetikzlibrary{automata}

\usepackage{amsfonts}
\usepackage{amssymb,latexsym,amsfonts,amsmath}
\usepackage{graphicx}

\newtheorem{theorem}{Theorem}[section]
\newtheorem{lemma}[theorem]{Lemma}

\newtheorem{proposition}[theorem]{Proposition}

\theoremstyle{definition}
\newtheorem{definition}[theorem]{Definition}
\newtheorem{example}[theorem]{Example}

\theoremstyle{remark}
\newtheorem{remark}[theorem]{Remark}
\numberwithin{equation}{section}

\renewcommand{\dim}{\mathrm{dim}}

\begin{document}
\begin{abstract}
We define observability and detectability for linear switching systems as the
possibility of reconstructing and respectively of asymptotically
reconstructing the hybrid state of the system from the knowledge of the
output for a suitable choice of the control input. We derive a necessary and
sufficient condition for observability that can be verified computationally. A
characterization of control inputs ensuring observability of switching
systems is given. Moreover, we prove that checking detectability of a linear switching
system is equivalent to checking asymptotic stability of a suitable switching
system with guards extracted from it, thus providing interesting links to
Kalman decomposition and the theory of stability of hybrid systems.
\end{abstract}

\title[Observability and Detectability of Linear Switching Systems: A Structural Approach]{Observability and Detectability of Linear Switching Systems: \\A Structural Approach}
\thanks{This work has been partially supported by the
HYCON Network of Excellence, contract number FP6-IST-511368 and by Ministero
dell'Istruzione, dell'Universita' e della Ricerca under Projects MACSI and
SCEF (PRIN05).}

\author[Elena De Santis, Maria Domenica Di Benedetto and Giordano Pola]{Elena De Santis, Maria Domenica Di Benedetto and Giordano Pola}
\address{Department of Electrical Information Engineering,Center of Excellence DEWS\\
University of L'Aquila, \\
Poggio di Roio, 67040 L'Aquila (Italy)}
\email{\{desantis,dibenede,pola\}@ing.univaq.it}
\maketitle

\section{Introduction}

Research in the area of hybrid systems addresses significant application domains with the aim of developing further understanding
of the implications of the hybrid model on control algorithms and to evaluate whether using this formalism can be of substantial
help in solving complex, real--life, control problems. In many application domains, hybrid controller synthesis problems are
addressed by assuming full hybrid state information, although in many realistic situations state measurements are not available.
Hence, to make hybrid controller synthesis relevant, the design of hybrid state observers is of fundamental importance. A step
towards a procedure for the synthesis of these observers is the analysis of observability and detectability of hybrid systems.
Observability has been extensively studied both in the continuous \cite{12-Kalman60,12-Luenberger71} and in the discrete domains
(see e.g. \cite{Ozv:AC90,12-Ramadge1986}). In particular, Sontag in \cite{12-SontagSIAM79} defined a number of observability
concepts and analyzed their relations for polynomial systems. More recently, various researchers investigated observability of
hybrid systems. The definitions of observability and the criteria to assess this property varied depending on the class of
systems under consideration and on the knowledge that is assumed at the output. Incremental observability was introduced in
\cite{12-bempAC00} for the class of piecewise affine systems. Incremental observability implies that different initial states
always give different outputs independently of the applied input.
 A characterization of observability and the definition of a
hybrid observer for the class of autonomous piecewise affine systems can be found in \cite{Collins04}. In \cite{Dinnocenzo:2006}
observability of autonomous hybrid systems was analyzed by using abstraction techniques. In \cite{12-Ball2003}, the notion of
generic final--state determinability proposed in \cite{12-SontagSIAM79} was extended to hybrid systems and sufficient conditions
were given for linear hybrid systems. The work in \cite{Vidal2003} considered autonomous switching systems and proposed a
definition of observability based on the concept of indistinguishability of continuous initial states and discrete state
evolutions from the outputs in free evolution. In \cite{DeS:CDC03,Babaali:2005} observability of switching systems (with control)
was investigated. Critical observability for safety critical switching systems was introduced in \cite{Des:LNCIS05}, where a set
of ``critical'' states must be reconstructed immediately since they correspond to hazards that may yield catastrophic events.\\
While observability of hybrid systems was addressed in the papers cited above, a general notion of detectability has not been
introduced as yet. To the best of our knowledge, the only contribution dealing with detectability can be found in
\cite{Mariton:dekker90} where detectability was defined for the class of jump linear systems as equivalent to the existence of a
set of linear
gains ensuring the convergence to zero of the estimation error in a stochastic setting.\\
In this paper we address observability and detectability for the class of switching systems.
General notions of observability and detectability are
introduced for the class of linear switching systems, though our definitions
apply to more general classes of hybrid systems, since they involve only
dynamical properties of the executions that are generated by the hybrid
system.
%The observability notion proposed here for linear switching systems
%and those proposed for other classes of hybrid systems in \cite{12-bempAC00},
%\cite{Vidal2003} are shown not to be equivalent in
%\cite{DeS:CDC03}.
Further, we derive a computable necessary and sufficient
condition for assessing observability.
%As a by--product, we obtain an explicit
%characterization of the class of control inputs ensuring observability and
%show that this class contains "almost all" inputs.
We then characterize detectability using a Kalman--like approach. In particular, we show that checking detectability of a linear
switching system is equivalent to checking asymptotic stability of a suitable linear switching system with guards associated with
the original system. This result is clearly related to the classical detectability analysis of linear systems. It is important
because it allows one to leverage a wealth of existing results on the stability of switching and hybrid systems (see e.g.
\cite{Narendra.tac1994,branicky:98,LiberzonBook} and the references therein). A preliminary version of this paper appeared before
in \cite{DeS:CDC03}. A characterization of observability, close to the one of \cite{DeS:CDC03} and of the one presented in this
paper, can be found in \cite{Babaali:2005} for a subclass of the switching systems considered in \cite{DeS:CDC03}.
The relation between \cite{Babaali:2005}, \cite{DeS:CDC03} and the present paper is discussed in Section 3. \\
The paper is organized as follows. In Section 2, we introduce linear switching systems and the notions of observability and
detectability. Section 3 is devoted to finding conditions for the reconstruction of the discrete component of the hybrid state.
In Section 4 we give a characterization of observability and detectability. In Section 5, an example shows the applicability and
the benefits of our results. Section 7 includes technical proofs of some of the results established in Section 3. Section 6
offers some concluding remarks.

\section{Preliminaries and basic definitions}

In this section, we introduce the notations and some basic definitions that are used in the paper.

\subsection{Notation}
The symbols $\mathbb{N}$,
$\mathbb{R}$ and $\mathbb{R}^{+}$ denote the natural, real and positive real
numbers, respectively. The symbol $I$ denotes the identity matrix of
appropriate dimensions. Given
a vector $x\in\mathbb{R}^{n}$, the symbol $x^{\prime}$ denotes the transpose
of $x$. The symbol $\Vert. \Vert_{n}$ denotes the Euclidean norm of a vector in the linear
space $\mathbb{R}^{n}$. Given a linear subspace $H$ of $\mathbb{R}^{n}$, the
symbol $dim(H)$ denotes its dimension and the symbol $\pi_{H}$ denotes the
projector on $H$, i.e. $\pi_{H}x$ is the Euclidean orthogonal
projection of $x$ onto $H$. Given a matrix $M\in\mathbb{R}^{n\times
m}$, the symbols $Im\left(  M\right)  $ and $ker\left(  M\right)  $ denote
respectively the range and the null space of $M$; given a set $H\subseteq\mathbb{R}^{n}$ the symbol $M^{-1}(H)$ denotes the inverse image of
$H$ through $M$, i.e. $M^{-1}(H)=\{x\in\mathbb{R}^{m}|\exists y\in H :
y=Mx\}$. Given a set $\Omega$, the symbol $card(\Omega)$ denotes the
cardinality of $\Omega$.

\subsection{Switching systems\label{sec:2.1}}

We consider the class of linear switching systems and the class of linear switching systems with guards, which generalize the
class defined in \cite{DeS:CDC03}, following the general model of hybrid automata of \cite{Lygeros:99,Tomlin:98}. Switching
systems are relevant in many application domains such as, among many others, mechanical systems, power train control, aircraft
and air traffic control, switching power converters, see e.g. \cite{LiberzonBook,Des:LNCIS05,Des:NA06} and the references
therein.
\newline The hybrid state $\xi$ of a
$GLSw$--system $\mathcal{H}$ is composed of two components: the discrete state
$i$ belonging to the finite set $Q=\{1,2,...,N\}$, called discrete state
space, and the continuous state $x$ belonging to the linear space
$\mathbb{R}^{n_{i}}$, whose dimension $n_{i}$ depends on $i\in Q$. The hybrid
state space of $\mathcal{H}$ is then defined by $\Xi=\bigcup_{i\in Q}\left\{
i\right\}  \times\mathbb{R}^{n_{i}}$. The control input of $\mathcal{H}$ is a
function $u\in\mathcal{U}$, where $\mathcal{U}$ denotes the class of piecewise
continuous functions $u:\mathbb{R}\rightarrow\mathbb{R}^{m}$. The output function of $\mathcal{H}$ belongs to the set $\mathcal{Y}$ of
piecewise continuous functions $y:\mathbb{R} \rightarrow \mathbb{R}^{l}$. The evolution of the continuous state $x$ and of the output $y$ of $\mathcal{H}$ is determined by the linear control systems:
\begin{equation}
S(i):\left\{
\begin{array}
[c]{l}%
\dot{x}=A_{i}x+B_{i}u,\\
y=C_{i}x,
\end{array}
\right.  \label{systema}%
\end{equation}
whose dynamical matrices $A_{i},B_{i},C_{i}$ depend on the current discrete state $i\in Q$. The evolution of the discrete state of
$\mathcal{H}$ is governed by a Finite State Machine (FSM), so that a
transition from a state $i\in Q$ to a state $h\in Q$ may occur if $e=(i,h)\in
E$, where $E\subseteq Q\times Q$ is the set of (admissible) transitions in the
FSM, and if the continuous state $x$ is in the set $G(e)\subseteq\mathbb{R}%
^{n_{i}}$, called guard\footnote{In this paper, the role of the guard $G(e)$ is to enable (and not to enforce) a transition.}
\cite{LygerosTAC03}. Whenever a transition $e=(i,h)$ occurs, the continuous state $x$ is instantly reset to a new value $R(e)x$,
where $R$ is the reset function which associates a matrix $R(e)\in\mathbb{R}^{n_{h}\times n_{i}}$ to each $e\in E$. We assume
that $R(e)\neq I$, for any in-loop transition
$e=(i,i)\in E$. \newline A \textit{linear switching system with guards} ($GLSw$%
--system) $\mathcal{H}$ is then specified by means of the tuple:
\begin{equation}
\left(  \Xi,S,E,G,R\right)  , \label{tuple}%
\end{equation}
with all the symbols as defined above. Given a \mbox{$GLSw$--system} $\mathcal{H}$, if $G(e)=\mathbb{R}^{n_{i}}$ for any $e\in E$, then
$\mathcal{H}$ is called \textit{linear switching system}\emph{\ }($LSw$--system) and for simplicity the symbol $G$ is omitted in
the tuple (\ref{tuple}), i.e. $\mathcal{H}=\left(  \Xi,S,E,R\right)  $. A $GLSw$--system $\mathcal{H}$ is said to be
\textit{autonomous} if all systems $S(i)$ are autonomous, i.e. \mbox{$B_{i}=0$}.
\newline The evolution in time of $GLSw$--systems can be defined as in
\cite{Lygeros:99}, by means of the notion of \textit{execution}. We recall that a hybrid time basis $\tau$ is an infinite or
finite sequence of sets $I_{j}=[t_{j},t_{j+1}),j=0,1,...,card(\tau)-1$, with $t_{j+1}>t_{j}$; let be $card(\tau)=L$. If
$L<\infty$, then $t_{L}=\infty$. Given a hybrid time basis $\tau$, time instants $t_{j}$ are called \textit{switching times}.
Throughout the paper we suppose that given a hybrid time basis, the number of switching times within any bounded time interval is
finite, thus avoiding Zeno behaviour \cite{LygerosTAC03} in the evolution of the system. Let $\mathcal{T}$ be the set of all
hybrid time bases and consider a collection:
\begin{equation}
\chi=\left(  \xi_{0},\tau,u,\xi,y\right)  , \label{exec}%
\end{equation}
where $\xi_{0}\in\Xi$ is the initial hybrid state, $\tau\in\mathcal{T}$ is the
hybrid time basis, $u\in\mathcal{U}$ is the continuous control input,
$\xi:\mathbb{R}\rightarrow\Xi$ is the hybrid state evolution and
$y\in\mathcal{Y}$ is the output evolution.
%The hybrid system temporal behaviour is defined by stating conditions such
%that $\chi$ represents a possible evolution of $\mathcal{H}$.
The function
$\xi$ is defined as follows:
%\[%
\[
\xi\left(  t_{0}\right)  =\xi_{0},\,\,\,\,\xi\left(  t\right)  =\left(
q(t),x(t)\right)  ,
\]
%\]
where at time $t\in I_{j}$, $q(t)=q(t_{j})$, $x(t)$ is the (unique) solution of the dynamical system $S(q(t_{j}))$, with initial
time $t_{j}$, initial state $x\left(  t_{j}\right)  $ and control law $u$. Moreover, if we set
${x^{-}}(t_{j})=lim_{t\rightarrow{t_{j}}^{-}}x(t)$ the following conditions have to be satisfied for any $j=1,...,L-1$:
\[%
\begin{array}
[c]{lll}%
(q(t_{j-1}),q(t_{j}))\in E, \\
{x^{-}}(t_{j})\in G(q(t_{j-1}),q(t_{j})), \\
x(t_{j})=R(q(t_{j-1}),q(t_{j})){x^{-}(}t_{j}).
\end{array}
\]
The output evolution $y$ is defined for any \mbox{$j=0,$}$1,...,$\mbox{$L-1$} by:%
\[
y(t)=C_{q(t_{j})}x(t),\hspace{0.5cm}t\in\lbrack t_{j},t_{j+1}).
\]
A tuple $\chi$ of the form (\ref{exec}), which satisfies the conditions above, is called an \textit{execution} of $\mathcal{H}$
\cite{LygerosTAC03}.

\subsection{Observability and Detectability\label{sec:2.2}}
%The definition of observability and detectability that we propose is based on
%the existence of at least an input--output experiment such that the hybrid state is reconstructed from some finite time on. \\
%In linear system theory, observability and detectability deal with the exact
%and (respectively) the asymptotic reconstruction of the state, on the basis of
%the knowledge of the continuous input and of the continuous output that is
%accessible from the environment.
In this section, we introduce the notions of observability and detectability for the class of \mbox{$GLSw-$systems}.\\
%Notice that in our framework, according to
%the definition of hybrid output, the transition from one discrete state to
%another may not be visible from the observed output.
%We start by introducing a metric on the hybrid state space of $GLSw$%
%--systems.
%For formally introducing the notions of observability and detectability,
Given a $GLSw$--system $\mathcal{H}$, we equip the hybrid state space with a metric:
% $\delta:\Xi\times\Xi\rightarrow\mathbb{R}^{+}\cup\{\infty\}$, such that, for any $(i,x_{i}),$ $(h,x_{h})\in\Xi$, $\delta((i,x_{i}),(h,x_{h})):=\left\Vert x_{i}-x_{h}\right\Vert _{n_{i}}$ if $n_{i}=n_{h}$ and $\delta((i,x_{i}),(h,x_{h})):=\infty$ if $n_{i}\neq n_{h}$.
\begin{equation}
\delta((i,x_{i}),(h,x_{h}))=\left\{
\begin{array}
[c]{ll}%
\infty, & \text{if }i\neq h,\\
\left\Vert x_{i}-x_{h}\right\Vert _{n_{i}}, & \text{if }i=h.
\end{array}
\right. \nonumber%
\end{equation}
The pair $(\Xi,\delta)$ is a metric space.
%Metric (\ref{delta}) has been obtained by combining classical metrics for discrete and Euclidean spaces.
\begin{definition}
\label{def_detect} A $GLSw$--system $\mathcal{H}$ is detectable if
there exist a control input $\widehat{u}\in\mathcal{U}$ and a function
$\widehat{\mathcal{\xi}}:\mathcal{Y}\times\mathcal{U}\rightarrow\Xi
$\textbf{\ }such that:%
\begin{equation}%
\begin{array}
[c]{l}%
\forall\varepsilon>0,\forall\rho>0,\exists\hat{t}>t_{0}:
\\
\delta(\widehat{\mathcal{\xi}}(\left.  y\right\vert _{\left[  t_{0},t\right]
},\left.  \widehat{u}\right\vert _{\left[  t_{0},t\right)  }),\xi\left(
t\right)  )\leq\varepsilon,
\\
\forall t\geq\hat{t},t\neq t_{j},j=0,1,...,L,
\end{array}
\label{cond_det}%
\end{equation}
for any execution $\chi$ with control input $\widehat{u}$ and hybrid initial
state $\xi_{0}=(i,x_{0})$ with $\left\Vert x_{0}\right\Vert _{n_{i}}\leq\rho
$. If condition (\ref{cond_det}) holds with $\varepsilon=0$, then
$\mathcal{H}$ is observable.
\end{definition}
By Definition \ref{def_detect}, an observable $GLSw$--system is also detectable.
%Definition\ \ref{def_detect} still holds for more general hybrid systems since
%it is based only on the notion of execution. A comparison between the
%observability notion of Definition \ref{def_detect} and some other
%observability notions available in the literature for other classes of hybrid
%systems (see e.g. \cite{12-bempAC00}, \cite{Vidal:cdc2002}, \cite{Vidal2003})
%can be found in \cite{DeS:CDC03}.
By specializing Definition \ref{def_detect} to linear systems, the classical observability and detectability notions are
recovered. Note that the reconstruction of the current hybrid state is required at every time $t\geq\hat{t}$ with $t\neq t_{j}$.
Time instants $t_{j}$ are ruled out as it is for observable linear systems, where the current state may be reconstructed only at
every time strictly greater than the initial time. However, observability and detectability for linear systems are defined
independently from the control function, while here we assume to choose a suitable control law. The two definitions coincide for
linear systems but not for $GLSw$--systems. In fact, if the observability (or detectability) property were required for
\emph{any} input function, then any $GLSw$--system would never be observable (or detectable), see e.g.
\cite{DeS:CDC03,Babaali:2005}.
%To see this, suppose
%that the initial continuous state is the origin and the control law is
%identically zero. Then, the output is identically zero and hence it is
%not possible to reconstruct the discrete state. This problem could be overcome by requiring the
%observability (or detectability) property to hold for any control input but
%not for all initial conditions in the state space, hence defining a
%\emph{local} notion of observability. We adopt a \textit{global} notion with respect to the initial states and we
However, we will show in Section 3 that if a
switching system is observable in the sense of Definition \ref{def_detect},
then it is observable for ``almost all'' input functions.
\newline
Definition \ref{def_detect} requires the reconstruction of the
discrete and of the continuous state. We consider these two issues separately, by stating conditions that ensure the
reconstruction of the discrete state in Section \ref{sec:dis_recon} and of the continuous state in Section \ref{sec:obs}.

\section{ Location observability\label{sec:dis_recon}}

In this section, we focus on the reconstruction of the discrete component of the hybrid state \textit{only}. By specializing
Definition \ref{def_detect}, we have:

\begin{definition}
\label{Def_loc_obs}A $GLSw$--system $\mathcal{H}$ is \textit{location observable} if there exist a control input $\hat{u}%
\in\mathcal{U}$ and a function $\widehat{q}:\mathcal{Y}\times\mathcal{U}%
\rightarrow Q$ such that:%
\begin{equation}%
\begin{array}
[c]{l}%
\forall\rho>0,\exists\hat{t}>t_{0}:
\\
\widehat{q}(\left.  y\right\vert _{\left[  t_{0},t\right]  },\left.
\widehat{u}\right\vert _{\left[  t_{0},t\right)  })=q\left(  t\right)  ,
\\
\forall t\geq\hat{t},t\neq t_{j},j=0,1,...,L-1,
\end{array}
\label{loc_obs_cond}%
\end{equation}
for any execution $\chi$ with control input $\widehat{u}$ and hybrid initial
state $\xi_{0}=(i,x_{0})$ with $\left\Vert x_{0}\right\Vert _{n_{i}}\leq
\rho$.
\end{definition}

A $GLSw$--system $\mathcal{H}$ is said to be \textit{location
observable} \textit{for a control input} $\hat{u}\in\mathcal{U}$ if there
exists a function \mbox{$\widehat{q}:\mathcal{Y}\times\mathcal{U}\rightarrow Q$} such
that condition (\ref{loc_obs_cond}) is satisfied.
%\newline By comparing
%Definitions \ref{def_detect} and \ref{Def_loc_obs}, it follows that location
%observability is a necessary condition for a $GLSw$--system to be observable
%or detectable.
The definition of location observability guarantees the
reconstruction of the discrete state, but not of the switching times,
as the following example shows.
%\footnote{This phenomenon will play an important role in the characterization of detectability, as discussed in the next section.}.
%Condition (\ref{condT}) of Lemma \ref{Th_nonempty} allows the reconstruction
%of the switching times, whenever a switching occurs between two different
%discrete states, but not in general if the transition resets the discrete
%state into itself, as the following example shows.
\begin{example}
\label{autociclo} Consider a $GLSw$--system $\mathcal{H}=(\Xi,S,$ $E,G,R)$, where
$\Xi=\left\{  1\right\}  \times\mathbb{R}^{3}$, $E=\{e\}$ with $e=(1,1)$ and
$G(e)=\mathbb{R}^{3}$. Let the dynamical system $S(1)$ and the reset function $R(e)$ be described by the
following dynamical matrices:
\[
A_{1}=\left(
\begin{array}
[c]{cc}%
1 & 0 \\
0 & 1 %
\end{array}
\right)  \text{,}B_{1}=\left(
\begin{array}
[c]{c}%
1\\
0
\end{array}
\right)  \text{,}C_{1}=\left(
\begin{array}
[c]{cc}%
1 & 0
\end{array}
\right)  \text{,}\\%
%\]
%with $\left(  A_{1},C_{1}\right)  $ detectable, $a_{22}>0$, $a_{12}\neq0$,
%and let the reset function be given by:
%\[
R(e)=\left(
\begin{array}
[c]{cc}%
1 & 0 \\
1 & 1
\end{array}
\right)  .
\]
The system $\mathcal{H}$ is trivially location observable for any control
input $u$. However since for any $x\in\mathbb{R}^{3}$, $\left(  R\left(
e\right)  -I\right)  x$ belongs to the kernel of the observability matrix
associated with $S(1)$, it is not possible to reconstruct the switching times,
for any choice of the control input $u$.
\end{example}

For later use, given \mbox{$i,h\in Q$}, define the following augmented linear system $S_{ih}$:%
\begin{equation}
\dot{z}=A_{ih}z+B_{ih}u\text{, \ \ \ \ }y_{ih}=C_{ih}z\text{,}
\label{sys_aumentato}%
\end{equation}
where:%
\[
\begin{array}
[c]{ccc}%
A_{ih}=\left(
\begin{array}
[c]{cc}%
A_{i} & 0\\
0 & A_{h}%
\end{array}
\right)  , & B_{ih}=\left(
\begin{array}
[c]{c}%
B_{i}\\
B_{h}%
\end{array}
\right)  , & C_{ih}=\left(
\begin{array}
[c]{cc}%
C_{i} & -C_{h}%
\end{array}
\right)  \text{.}%
\end{array}
\]
Let $\mathcal{V}_{ih}\subseteq\mathbb{R}^{n_{i}+n_{h}}$ be the maximal
controlled invariant subspace \cite{Basile} for system $S_{ih}$ contained in
$\ker(C_{ih})$, i.e. the maximal subspace \mbox{$F\subseteq\mathbb{R}^{n_{i}+n_{h}}$}
satisfying the following sets inclusions:%
\begin{equation}%
\begin{array}
[c]{ccc}%
A_{ih}F\subseteq F+Im(B_{ih}), &  & F\subseteq ker (C_{ih}).
\end{array}
\label{inv}%
\end{equation}
Define $\hat{J}=\left\{  \left(  i,h\right)\in Q\times Q  : i\neq h\right\}  $
and
%\[
%\begin{array}{l}
%J_{p}=\{i\in Q:\gamma _{Q}(i)=p\},p\in P, \\
%\hat{J}=\left\{ \left( i,h\right) :i,h\in J_{p},i\neq h,p\in P\right\} ,%
%\end{array}%
%\]%
%and
consider the set:
\[
\mathcal{U}^{\ast}=\left\{  u\in\mathcal{U}:u\neq\widetilde{u},\text{
}a.e.,\text{ }\forall\widetilde{u}\in\widetilde{\mathcal{U}}\right\}  ,
\]
where:%
\begin{equation}%
\begin{array}
[c]{l}%
\widetilde{\mathcal{U}}=\bigcup_{\left(  i,h\right)  \in\hat{J}}%
\mathcal{U}_{ih},\\
\mathcal{U}_{ih}=\left\{
\begin{array}
[c]{cc}%
u\in\mathcal{U}:u(t)=K_{ih}z(t)+v_{ih}%
(t),\\
t\geq\hat{t}\text{, for some }\hat{t}\in\mathbb{R}
\end{array}
\right\}  ,
\end{array}
\label{Z2}%
\end{equation}
the gain $K_{ih}$ is such that $\left(  A_{ih}+B_{ih}K_{ih}\right)
\mathcal{V}_{ih}\subseteq\mathcal{V}_{ih}$, $v_{ih}(t)\in B_{ih}^{-1}\left(
\mathcal{V}_{ih}\right)  ,\forall t\geq\hat{t}$ and $z(t)$ is the state
of system $S_{ih}$ at time $t$, under control $u$ with $z(\hat{t}%
)\in\mathcal{V}_{ih}$.
%starting from some $z(\overline{t}%
%)\in\mathcal{V}_{ih}$, under control $\left.  u\right\vert _{\left[
%\overline{t},t\right)  }$.
%The set \ $\mathcal{U}_{ih}$ is composed
%by those control inputs $u\in\mathcal{U}$ so that after a finite time
%$\hat{t}$ the corresponding output $y_{ih}$ of $S_{ih}$ with
%$z(\hat{t})\in\mathcal{V}_{ih}$ is identically zero.
%ure
%\textit{zero dynamics} \cite{Kailath:80}\ to system $S_{ih}$.
%If $B_{ih}%
%^{-1}\left(  \mathcal{V}_{ih}\right)  =\left\{  0\right\}  $, i.e. if the
%system $S_{ih}$ is left invertible, then
%\[
%\mathcal{U}_{ih}=\left\{  \left.  u\right\vert _{\left[  \overline{t}%
%,\infty\right)  }:\overline{t}\in\mathbb{R}\text{, }u(t)=\exp\left(
%A_{ih}+B_{ih}K_{ih}\right)  \left(  t-\overline{t}\right)  z(\overline
%{t}),\text{ }z(\overline{t})\in\mathcal{V}_{ih},t\geq\overline{t}\right\}
%\]
%Therefore, if $\forall\left(  i,h\right)  \in\hat{J}\mathbf{\ }$the system
%$S_{ih}$ is left invertible and $\mathbf{\exists}k\in\mathbb{N},k<n_{i}%
%+n_{h}:C_{i}A_{i}^{k}B_{i}\neq C_{h}A_{h}^{k}B_{h}$, then obviously the set
%$\mathcal{U}^{\ast}$ is nonempty. The next theorem states that the above
%condition on the input-output behavior of the systems $S\left(  i\right)  $
%and $S\left(  h\right)  $ is actually sufficient for nonemptyness of
%$\mathcal{U}^{\ast}$ in the general case.
%The following result shows that, under appropriate conditions on the system
%parameters, the set $\mathcal{U}^{\ast}$ is nonempty.
%Before giving the main result of this section we need the following technical result.
The set $\mathcal{U}^{*}$ is composed of the control inputs $u$ such that after a finite time $\hat{t}$ the output $y_{ih}$ of
$S_{ih}$ with any initial state \mbox{$x_{0}\in\mathbb{R}^{n_{i}+n_{h}}$} and the control input $u$ is not identically zero for any
choice of $(i,h)\in \hat{J}$. We will show that control inputs in $\mathcal{U}^{*}$ ensure the reconstruction of the discrete
state. The following result identifies conditions for nonemptyness of $\mathcal{U}^{*}$.

\begin{lemma}
\label{Th_nonempty}Given a $GLSw-$system $\mathcal{H}$, the set $\mathcal{U}%
^{\ast}$ is nonempty if
\begin{equation}
\forall\left(  i,h\right)  \in\hat{J}\mathbf{,\exists}k\in\mathbb{N}%
,k<n_{i}+n_{h}:C_{i}A_{i}^{k}B_{i}\neq C_{h}A_{h}^{k}B_{h}. \label{condT}%
\end{equation}

\end{lemma}

%\begin{proof}
%see Appendix.
%\end{proof}
%
%\bigskip
The proof of the above result requires some technicalities and is therefore reported in the Appendix. We now have all the
ingredients for characterizing location observability of switching systems.

\begin{theorem}
\label{Th_locobs} A $GLSw$--system $\mathcal{H}$ is location observable if and only if
condition (\ref{condT}) holds.
\end{theorem}

\begin{proof}
(Necessity) Suppose by contradiction, that $\exists\left(  i,h\right)
\in\hat{J}$ such that condition (\ref{condT}) is not satisfied and consider any $u\in\mathcal{U}$ and any executions $\chi_{1}=\left(  (i,0),\tau,u,\xi_{1},y_{1}\right)$ and $\chi_{2}=\left(
(h,0),\tau,u,\xi_{2},y_{2}\right)$ with $\tau=\{I_{0}\}$ and $I_{0}=[0,\infty)$. It is readily seen that
$y_{1}=y_{2}$ and therefore the discrete state cannot be
reconstructed. (Sufficiency)
%If $card(J_{p})\leq 1$, $\forall
%p\in P$, then the current discrete state can be reconstructed, on the basis
%of the discrete component of the observed output. Therefore assume $\exists
%p\in P$ such that $card(J_{p})>1$.
By Lemma \ref{Th_nonempty}, condition (\ref{condT}) implies that
$\mathcal{U}^{\ast}\neq\varnothing$; choose any $u\in\mathcal{U}^{\ast}$ and consider any execution $\chi=(\xi_{0},\tau,u,\xi,y)$. Consider any $j<L$ and
let $\xi(t)=(i,x(t)),t\in\lbrack t_{j},t_{j+1})$. Given any $h\in Q$, denote
by $y_{ih}(t,t_{j},z,\left.  u\right\vert _{\left[  t_{j},t\right)  })$ the
output evolution at time $t$ of system $S_{ih}$ with initial state
$z\in\mathbb{R}^{n_{i}+n_{h}}$ at initial time $t_{j}$ and control law
$\left.  u\right\vert _{\left[  t_{j},t\right)  }$. Since $u\in\mathcal{U}%
^{\ast}$ then for any $\varepsilon>0$, for any $h\neq i$ and for any
$w\in\mathbb{R}^{n_{h}}$ there exists a time $t\in(t_{j},t_{j}+\varepsilon)$
such that $y_{ih} (  t,t_{j},(
\begin{array}
[c]{ll}%
x(t_{j}) &
w
\end{array}
)'  ,u
%\right\vert _{\left[  t_{j},t\right)}
)  \neq0
$.
%By definition of $\mathcal{U}_{ih}$ as in (\ref{Z2}), $\exists\varepsilon>0$,
%such that $y_{ih}(t,t_{j},x_{t_{j}},\left.  u\right\vert _{\left[
%t_{j},t\right)  })\neq0$, for any $t\in\left(  t_{j+1},t_{j+1}+\varepsilon
%\right]  $, for any $x_{t_{j}}$, and for any $h\neq i$.
This implies that $y(t)\neq y_{h}(t)$, where $y_{h}$ is the output associated with the execution
$(\xi_{0h},\tau,u,\xi_{h},y_h)$ with $\xi_{h}(t)=(h,x_{h}(t)),t\in [t_{j},t_{j+1})$. Hence, the discrete state can be reconstructed for any
\mbox{$t\in(t_{j},t_{j+1})$}, and the statement follows.
\end{proof}
It is seen from the above result that \textit{if a }$GLSw$--system $\mathcal{H}$\textit{ is location observable then it is
location observable for any input function }$u\in\mathcal{U}^{\ast}$. A control law that ensures location observability is
derived in the proof of Lemma \ref{Th_nonempty}. Moreover, if the set of control inputs is the set
$\mathcal{C}^{\infty}(\mathbb{R}^{m})$ of smooth functions $u:\mathbb{R}\rightarrow \mathbb{R}^{m}$ (instead of the set
$\mathcal{U}$ of piecewise continuous functions), then $\mathcal{U}^{\ast}$ contains all and nothing but the control inputs which
ensure location observability.
%The following result shows that any control input in $\mathcal{U}^{*}$ ensures the reconstruction of the discrete state. The proof is a straightforward consequence of the proof of Theorem \ref{Th_locobs} and is therefore omitted.
%\begin{corollary}
%If a $GLSw$--system $\mathcal{H}$ is location observable then it is
%location observable for any input function $u\in\mathcal{U}^{\ast}$.
%\end{corollary}
\begin{remark}
Condition (\ref{condT}) was first given in \cite{DeS:CDC03} as a necessary and sufficient condition for guaranteeing location
observability of linear switching systems. A subclass of switching systems was then considered in \cite{Babaali:2005} where
similar observability conditions can be found. While the notion of observability of \cite{Babaali:2005} and the one in the
present paper (Definition \ref{def_detect} or equivalently the definition in \cite{DeS:CDC03}) are slightly different, the
notions of location observability coincide in the two papers. This translates in a characterization of location observability in
\cite{Babaali:2005} which is equivalent to the one in \cite{DeS:CDC03} and hence to the one of the present paper (compare Theorem
3 of \cite{Babaali:2005}, Theorem 8 of \cite{DeS:CDC03} and Theorem \ref{Th_locobs} of this paper).
\end{remark}

\section{Characterizing Observability and Detectability\label{sec:obs}}
Definition \ref{def_detect} implies that \textit{a }$GLSw$\textit{--system is observable if and only if it is location observable
and }$S\left(  i\right)  $\textit{ is observable for any }$i\in Q$. \\
The intuitive algorithm for the reconstruction of the (current) hybrid state of an observable $GLSw$--system $\mathcal{H}$,
processes the output $y\in\mathcal{Y}$ and the input $u\in\mathcal{U}^{*}$. It first reconstructs the current discrete state, by
looking for the unique $i\in Q$ such that\footnote{If the switching system $\mathcal{H}$ is location observable and \mbox{$u\in
\mathcal{U}^{*}$}, Theorem \ref{Th_locobs} guarantees that such discrete state $i$ is unique.}
\begin{equation}
Y^{(n_{i})}(t)\in Im(\mathcal{O}_{i})+\mathcal{F}_{i}u(t),
\label{DiscreteObs}
\end{equation}
where $Y^{(n_{i})}(t)=(
\begin{array}
[c]{cccc}%
y(t)' & \dot{y}(t)' & \ldots & y^{(n_{i}-1)}(t)'
\end{array}
)'$, $\mathcal{O}_{i}$ is the observability matrix associated with $S(i)$ and
\begin{equation}
\mathcal{F}_{i}=\left(
\begin{array}
[c]{cccc}%
C_{i} & 0 & \ldots & 0\\
C_{i}A_{i} & C_{i}B_{i} & \ldots & 0\\
\ldots & \ldots & \ldots & 0\\
C_{i}A_{i}^{n_{i}} & C_{i}A_{i}^{n_{i}-1}B & \ldots & C_{i}B_{i}
\end{array}
\right);\nonumber
\end{equation}
Then, on the basis of the knowledge of $i$, it reconstructs the current continuous state $x(t)$, by computing:
\begin{equation}
\{x(t)\}=\mathcal{O}_{i}^{-1}\left(Y^{(n_{i})}(t)-\mathcal{F}_{i}u(t)\right).
%X^{(n)}(t)=\mathcal{O}_{i}^{-1}\left(Y^{(n)}(t)\right),
\label{ContinuousObs}
\end{equation}

%\newline
We now focus on $LSw$--systems and derive conditions that ensure detectability.
Since location observability is a necessary condition for a switching system to be observable or detectable, we assume now that
this property holds for all systems considered in this section. Given a $LSw$--system $\mathcal{H}=\left(  \Xi,S,E,R\right)$,
define the autonomous $LSw$--system:
\begin{equation}
\mathcal{H}^{\prime}=\left(  \Xi,S^{\prime},E,R\right)  , \label{S_primo}%
\end{equation}
where $S^{\prime}(i)$ is defined as $S(i)$ in (\ref{systema}) with $B_{i}=0$. We assume that $\mathcal{H}^{\prime}$ is with full
discrete evolution information, i.e. that the discrete state and the switching times are known at any time. Clearly,
detectability of $\mathcal{H}$ implies detectability of $\mathcal{H}^{\prime}$. Under some appropriate conditions, the converse
implication is true:

\begin{lemma}
\label{Prop2}A location observable $LSw$--system $\mathcal{H}$ is detectable if
$\mathcal{H}^{\prime}$ is detectable and $\mathcal{H}$ satisfies the following
property:%
\begin{equation}
E^{\circlearrowright}=\varnothing\text{ \ \ \ \ or \ \ \ \ }Im(R(e)-I)\cap
\ker(\mathcal{O}_{i})=\{0\},\forall e\in E^{\circlearrowright}, \label{Co}%
\end{equation}
where $E^{\circlearrowright}=\{(i,h)\in E: i=h\}$ and
$\mathcal{O}_{i}$ is the observability matrix associated with $S(i)$.
%with the linear system $S(i)$
%For any $i\in Q$, let $\mathcal{O}_{i}$ be the observability matrix associated
%with the linear system $S(i)$.
%then $\mathcal{H}$ is detectable.
\end{lemma}

Under condition (\ref{Co}), if a transition $(i,i)\in E^{\circlearrowright}$ occurs in $\mathcal{H}$ at time $t_{j}$ from a
hybrid state $(i,x^{-})$ to a hybrid state $(i,x^{+})$ with $x^{+}=R(i,i)x^{-}\neq x^{-}$ then $x^{+}-x^{-}\notin
ker(\mathcal{O}_{i})$. Hence the switching time $t_{j}$ can be reconstructed\footnote{Note that the switching system of Example
\ref{autociclo} does not satisfy condition (\ref{Co}) and therefore switching times in that case cannot be reconstructed.}. Then,
the proof of the result above just follows from the linearity of the continuous dynamics in $\mathcal{H}$ and from the definition
of $\mathcal{H}^{\prime}$.\\
The result of Lemma \ref{Prop2} reduces the analysis of detectability of a linear switching system \textit{with control}, to that
of an \textit{autonomous} linear switching system.\\
For analyzing detectability of $\mathcal{H}^{\prime}$ it is useful to first
perform a discrete state space decomposition. \\
Given $\mathcal{H}^{\prime
}=\left(  \Xi,S^{\prime},E,R\right)  $ as in (\ref{S_primo}) and a set
$\hat{Q}\subseteq Q$ let
\[
\left.  \mathcal{H}^{\prime}\right\vert _{\hat{Q}}=(\left.  \Xi\right\vert
_{\hat{Q}},\left.  S^{\prime}\right\vert _{\hat{Q}}\left.  ,E\right\vert
_{\hat{Q}},
\left.
R\right\vert_{\hat{Q}}),
\]
be the switching sub--system of $\mathcal{H}^{\prime}$ obtained by restricting
the discrete state space $Q$ of $\mathcal{H}$ to $\hat{Q}$, i.e. such that
$\left.  \Xi\right\vert _{\hat{Q}}=\bigcup_{i\in\hat{Q}}\left\{  i\right\}
\times\mathbb{R}^{n_{i}}$, $\left.  S^{\prime}\right\vert _{\hat{Q}}(i)=S^{\prime}(i),%\forall i\in\hat{Q}
$
$\left.  E\right\vert _{\hat{Q}}=\{(i,h)\in E:i,h\in\hat{Q}\}$ and $\left.
R\right\vert _{\hat{Q}}(i,h)=R(i,h)%\forall(i,h)\in\left.  E\right\vert _{\hat{Q}}
$.
%It shows that the system $\mathcal{S}^{\prime}$ is detectable if and
%only if a linear switching system, obtained by skipping all discrete states
%associated with observable linear systems of $\mathcal{S}^{\prime}$, is detectable.
\begin{proposition}
\label{Th_dec_dis2}The $LSw$--system $\mathcal{H}^{\prime}$ is detectable if and only if the $LSw$--system
$\left.  \mathcal{H}^{\prime}\right\vert _{\widehat{Q}}$ with $\widehat
{Q}=\{i\in Q:S(i)$ is not observable$\}$ is detectable.
\end{proposition}

\begin{proof}
(Necessity) Obvious. (Sufficiency) Consider any execution $\chi$ of $\mathcal{H}^{\prime}$. If $q(t)\in\widehat{Q}$ for any time
$t\geq t_{0}$ then the detectability of $\left.  \mathcal{H}^{\prime }\right\vert _{\widehat{Q}}$ implies the asymptotic
reconstruction of the hybrid state evolution of $\chi$. If $q(t)\notin\widehat{Q}$ for some finite time $t$, then
$S^{\prime}(q(t))$ is observable and hence it is possible to (exactly) reconstruct the continuous state of $\mathcal{H}^{\prime}$
in infinitesimal time. Once the continuous state $x(t')$ is known  at time $t'>t$, location observability of
$\mathcal{H}^{\prime}$ ensures the reconstruction of the hybrid state for any time $t''\geq t'$ with $t''\neq t_{j}$.
\end{proof}

%\newline
By Proposition \ref{Th_dec_dis2} there is no loss of generality in assuming that system $S^{\prime}(i)$ is not observable for any
$i\in Q$. Moreover, we assume that $S^{\prime}(i)$, $i\in Q$, are in observability canonical form, i.e. that dynamical matrices
associated with $S^{\prime }(i)$\ are of the form:
\[
A_{i}=\left(
\begin{array}
[c]{cc}%
A_{i}^{(11)} & 0\\
A_{i}^{(21)} & A_{i}^{(22)}%
\end{array}
\right)  ,C_{i}=\left(
\begin{array}
[c]{cc}%
C_{i}^{(1)} & 0
\end{array}
\right)  ,
\]
where $A_{i}^{(22)}\in\mathbb{R}^{d_{i}\times d_{i}}$, $0<d_{i}\leq n_{i}$ matrices $A_{i}^{(11)}$, $A_{i}^{(21)}$\ are of
appropriate dimensions and $(A_{i}^{(11)},C_{i}^{(1)})$ is an observable matrix pair, for any $i\in Q$. This assumption is made
without loss of generality: suppose that, for some $i\in Q$, the dynamical matrices $A_{i},C_{i}$ of the switching system
$\mathcal{H}^{\prime}$ are not in the observability canonical form. Then, we define an invertible linear transformation
$T_{i}:\mathbb{R}^{n_{i}}\rightarrow\mathbb{R}^{n_{i}} $ such that $T_{i}A_{i}T_{i}^{-1}$ and $C_{i}T_{i}^{-1}$ are in the
observability canonical form. For all $j\in Q$ such that the dynamical matrices $A_{j},C_{j}$ of the switching system
$\mathcal{H}^{\prime}$ are in the observability canonical form, we let $T_{j}$ be the identity matrix. We then define the
\textit{hybrid state space transformation} $T:\Xi\rightarrow\Xi$ such that for any $(i,x)\in\Xi$, $T(i,x):=(i,T_{i}x)$. The reset
function in the new coordinates is given by $T_{h}R(e)T_{i}^{-1}$, for any $e=(i,h)\in E$. 
The continuous component $x$ of the hybrid state $(i,x)$ of $\mathcal{H}%
^{\prime}$ can be partitioned as $x=(%
\begin{array}
[c]{cc}%
x_{1}^{\prime} & x_{2}^{\prime}%
\end{array}
)^{\prime}$, with $x_{1}\in\mathbb{R}^{n_{i}-d_{i}}$, $x_{2}\in
\mathbb{R}^{d_{i}}$, and the reset matrix $R(e)$ can be
partitioned as:%
\[
R(e)=\left(
\begin{array}
[c]{cc}%
R^{(11)}(e) & R^{(12)}(e)\\
R^{(21)}(e) & R^{(22)}(e)
\end{array}
\right)  ,
\]
where $R^{(22)}(e)\in\mathbb{R}^{d_{h}\times d_{i}}$ and $R^{(11)}(e)$, $R^{(12)}(e)$, $R^{(21)}(e)$\ are of appropriate
dimensions. Given the $LSw$--system $\mathcal{H}^{\prime}$ as in (\ref{S_primo}), define the $GLSw$--system:
\begin{equation}
\mathcal{H}_{0}=\left(  \Xi_{0},S_{0},E,G_{0},R_{0}\right)  , \label{S0}%
\end{equation}
where:

\begin{itemize}
\item $\Xi_{0}=\bigcup_{i\in Q}\left\{  i\right\}  \times\mathbb{R}^{d_{i}}$;

\item $S_{0}(i)$ is described by dynamics $\dot{z}(t)=A_{i}^{(22)}z(t)$, for
any $i\in Q$;

\item $G_{0}(e)=\ker(R^{(12)}(e))$, for any $e\in E$;

\item $R_{0}(e)=R^{(22)}(e)$, for any $e\in E$.
\end{itemize}

There is a strong connection between detectability of $\mathcal{H}^{\prime}$ and asymptotic stability of $\mathcal{H}_{0}$.
%We first give a
%characterization of detectability of $\mathcal{H}^{\prime}$ and then we
%establish this connection between stability of $\mathcal{H}_{0}$ and
%detectability of $\mathcal{H}^{\prime}$.
%as defined in (\ref{S_primo}).
%An execution $\chi=\left(  \xi_{0},\tau,\mathbf{0},\xi,y\right)  $ of
%$\mathcal{H}^{\prime}$ is said to be a \textit{hidden execution} if the output function $y$ is identically zero, i.e. if:
%\[
%C_{q(t)}x(t)=0,\forall t\geq t_{0}.
%\]
%%where for any $j=0,1,...,L-1$, $C_{i}$ is the output matrix associated with
%%the current discrete state $q(j)=i$.
Set \mbox{$\mathcal{B}:=\bigcup\nolimits_{i\in Q}\left\{  i\right\}
\times\mathcal{B}_{i}$}, where $\mathcal{B}_{i}=\{x\in\mathbb{R}^{n_{i}%
}:\left\Vert x\right\Vert _{n_{i}}\leq1\}$. We also define $\varepsilon\mathcal{B}:=\bigcup\nolimits_{i\in Q}\left\{  i\right\}
\times\varepsilon\mathcal{B}_{i}$ for any $\varepsilon\in \mathbb{R}^{+}$. 
An autonomous $GLSw$--system $\mathcal{H}$ is \textit{asymptotically stable} if the
continuous component of the hybrid state of any execution $\chi$ of
$\mathcal{H}$ converges to the origin as time goes to infinity, or
equivalently:
\[
\forall\varepsilon>0,\forall\rho>0,\exists\hat{t}>t_{0}:\xi(t)\in
\varepsilon\mathcal{B},\forall t\geq\hat{t},
\]
for any execution $\chi$ with hybrid initial state $\xi_{0}\in\rho\mathcal{B}%
$.
%\end{definition}
The following holds:

\begin{proposition}
\label{Prop_CNSdetect}The $LSw$--system $\mathcal{H}^{\prime}$ is detectable if
and only if the $GLSw$--system $\mathcal{H}_{0}$ is asymptotically stable.
\end{proposition}

\begin{proof}
(Sketch.) Let $\mathcal{E}_{0}$ be the set of executions of $\mathcal{H}^{\prime}$ such that $C_{q(t)}x(t)=0,\forall t\geq
t_{0}$. The continuous component $x(t)$ of the hybrid state $(q(t),x(t))$ of any execution in $\mathcal{E}_{0}$ belongs to the
subspace $\ker\left(  \mathcal{O}_{i}\right)  $ with $i=q(t)$ for any $t\in I_{j}$ and \mbox{$j=0,1,...,L$}. By definition of
$\mathcal{E}_{0}$, $\mathcal{H}^{\prime}$ is detectable if and only if the continuous component of the hybrid state $\xi$ of any
$\chi\in\mathcal{E}_{0}$ converges to the origin, i.e. $\forall\varepsilon>0,\forall\rho>0,\exists\hat{t}\geq t_{0}$ such that
$\xi(t)\in \varepsilon\mathcal{B},\forall t\geq\hat{t}$, for any $\chi\in\mathcal{E}_{0}$ with hybrid initial state
$\xi_{0}\in\rho\mathcal{B}$. By definition of the observability canonical form, this is equivalent to asymptotic stability of
$\mathcal{H}_{0}$.
\end{proof}

%\newline
By combining Lemma \ref{Prop2} and Propositions \ref{Th_dec_dis2} and \ref{Prop_CNSdetect} we obtain the following
characterization of detectability of $LSw$--systems.

\begin{theorem}
\label{Th_all}A $LSw$--system $\mathcal{H}$ is detectable if the following
conditions are satisfied:

\begin{description}
\item[i)] $\mathcal{H}$ is location observable, and

\item[ii)] $\mathcal{H}$ satisfies condition (\ref{Co}), and

\item[iii)] $\mathcal{H}_{0}$ is asymptotically stable.
\end{description}

Conversely, if $\mathcal{H}$ is detectable then conditions i) and iii) are
satisfied. \label{thdetect}
\end{theorem}

Since the executions associated with a $GLSw$--system $\left(\Xi,S,E,G,R\right)$ are also executions of the $LSw$--system
$\mathcal{H}=\left(\Xi,S,E,R\right)$, the conditions of Theorem \ref{Th_all} are also sufficient for a $GLSw$--system to be
detectable.

Detectability of switching systems has also been addressed in \cite{DeS:CDC03}. The above result provides a deeper analysis than
the one in \cite{DeS:CDC03} since it reduces detectability of $LSw$--systems to asymptotic stability of $GLSw$--systems (compare
Theorem 9 of \cite{DeS:CDC03} with the above result). This allows one to leverage the rich literature on stability of hybrid
systems (see e.g. \cite{Narendra.tac1994,branicky:98,LiberzonBook} and the references therein) for checking detectability. While
checking conditions i) and ii) is straightforward, checking condition iii) requires the analysis of asymptotic stability of
switching systems with guards.\\
We now derive sufficient conditions for assessing the asymptotic stability of $\mathcal{H}_{0}$, by \textit{abstracting}
$\mathcal{H}_{0}$ with linear switching systems with no guards. Given the autonomous $GLSw$--system $\mathcal{H}_{0}$ as in
(\ref{S0}) define the following autonomous $LSw$--systems:
\begin{equation}
\mathcal{H}_{1}=\left(  \Xi_{0},S_{0},E,R_{0}\right)  ,\hspace{0.5cm}
\mathcal{H}_{2}=\left(  \Xi_{0},S_{0},E,R_{2}\right)  ,
\label{abs}
\end{equation}
where $R_{2}(e)=R^{(22)}(e)\pi_{\ker(R^{(12)}(e))}$. The following result holds:

\begin{proposition}
\label{Th_detect}The autonomous $GLSw$--system $\mathcal{H}_{0}$ is asymptotically stable if either $\mathcal{H}_{1}$ or
$\mathcal{H}_{2}$ is asymptotically stable.
\end{proposition}
Since transitions in $LSw$--systems $\mathcal{H}_{1}$ and $\mathcal{H}_{2}$ are independent of the continuous state, the
asymptotic stability analysis of $\mathcal{H}_{1}$ and $\mathcal{H}_{2}$ is in general easier than the one of $\mathcal{H}_{0}$.
%The following examples show an application of the result above and point out
%that asymptotic stability of $\mathcal{S}_{1}\ $does not imply asymptotic
%stability of $\mathcal{S}_{2}$ and vice versa.
An application of this result is shown in the next section.

\section{An illustrative example}

In this section, we present an example that shows the interest and applicability of our results. 
Consider the linear switching system $\mathcal{H}=\left( \Xi,S,E,R\right)$, where:

\begin{itemize}
\item $\Xi=\left(  \{1\}\times\mathbb{R}^{4}\right)  \cup\left(
\{2\}\times\mathbb{R}^{3}\right)  \cup\left(  \{3\}\times\mathbb{R}%
^{2}\right)  \cup\left(  \{4\}\times\mathbb{R}\right) \cup\left(  \{5\}\times\mathbb{R}^{3}\right) \cup\left(  \{6\}\times\mathbb{R}^{2}\right) $;

\item $S$ associates to any $i\in Q=\{1,2,3,4,5,6\}$ the linear control
system $S(i)$ of (\ref{systema}), where:
%\[
%S(i):\left\{
%\begin{array}
%[c]{l}%
%\dot{x}(t)=A_{i}x(t)+B_{i}u(t),\\
%y(t)=C_{i}x(t),
%\end{array}
%\right.
%\]
%where:

\begin{equation}
\begin{array}
[c]{llll}
A_{1}=\left(
\begin{array}
[c]{rrrr}%
1 & 2 & 0 & 0\\
0 & 1 & 0 & 0\\
0 & 0 & -2 & 1\\
0 & 0 & 1 & -2
\end{array}
\right)  , &  B_{1}=\left(
\begin{array}
[c]{c}%
1\\
0\\
2\\
1
\end{array}
\right)  , & C_{1}=\left(
\begin{array}
[c]{cccc}%
1 & 1 & 0 & 0
\end{array}
\right)  ,\nonumber\\
A_{2}=\left(
\begin{array}
[c]{rrr}%
2 & 0 & 0\\
0 & -1 & 1\\
0 & 1 & -2
\end{array}
\right)  , &  B_{2}=\left(
\begin{array}
[c]{c}%
1\\
1\\
1
\end{array}
\right)  , & C_{2}=\left(
\begin{array}
[c]{ccc}%
1 & 0 & 0
\end{array}
\right)  ,\nonumber\\
A_{3}=\left(
\begin{array}
[c]{rrrr}%
1 & 0  \\
1 & -1
\end{array}
\right)  , &  B_{3}=\left(
\begin{array}
[c]{c}%
0\\
0
\end{array}
\right), & C_{3}=\left(
\begin{array}
[c]{cc}%
1 & 0
\end{array}
\right)  ,\nonumber\\
A_{4}=3, &  B_{4}=1, & C_{4}=1,\nonumber\\
A_{5}=\left(
\begin{array}
[c]{rrr}%
1 & 0 & 0  \\
1 & -1 & 0 \\
1 & 0 & -2
\end{array}
\right)  , &  B_{5}=\left(
\begin{array}
[c]{c}%
4\\
0\\
0
\end{array}
\right), & C_{5}=\left(
\begin{array}
[c]{ccc}%
1 & 0 & 0
\end{array}
\right)  ,\nonumber\\
A_{6}=\left(
\begin{array}
[c]{rrrr}%
5 & 0  \\
2 & -3
\end{array}
\right)  , &  B_{6}=\left(
\begin{array}
[c]{c}%
1\\
0
\end{array}
\right), & C_{6}=\left(
\begin{array}
[c]{cc}%
1 & 0
\end{array}
\right)  ;\nonumber\\

\end{array}
\end{equation}

%\begin {eqnarray}%
%&
%%\begin{align}
%A_{1}=\left(
%\begin{array}
%[c]{rrrr}%
%1 & 2 & 0 & 0\\
%0 & 1 & 0 & 0\\
%0 & 0 & -2 & 1\\
%0 & 0 & 1 & -2
%\end{array}
%\right)  , &  B_{1}=\left(
%\begin{array}
%[c]{c}%
%1\\
%0\\
%2\\
%1
%\end{array}
%\right)  ,C_{1}=\left(
%\begin{array}
%[c]{cccc}%
%1 & 1 & 0 & 0
%\end{array}
%\right)  ,\nonumber\\
%&
%A_{2}=\left(
%\begin{array}
%[c]{rrr}%
%2 & 0 & 0\\
%0 & -1 & 1\\
%0 & 1 & -2
%\end{array}
%\right)  , &  B_{2}=\left(
%\begin{array}
%[c]{c}%
%1\\
%1\\
%1
%\end{array}
%\right)  ,C_{2}=\left(
%\begin{array}
%[c]{ccc}%
%1 & 0 & 0
%\end{array}
%\right)  ,\nonumber\\
%&
%A_{3}=\left(
%\begin{array}
%[c]{rrrr}%
%1 & 0  \\
%1 & -1
%\end{array}
%\right)  , &  B_{3}=\left(
%\begin{array}
%[c]{c}%
%0\\
%0
%\end{array}
%\right)  ,C_{3}=\left(
%\begin{array}
%[c]{cc}%
%1 & 0
%\end{array}
%\right)  ,\nonumber\\
%&
%A_{4}=3, &  B_{4}=1,C_{4}=1;\nonumber
%%\end{align}
%\end{eqnarray}

\item $E=\{(1,2),(2,1),(2,3),(2,5),(3,3),(3,6),(4,1),(4,2),$ $(5,4),(5,6),(6,5)\}$;

\item $R$ is defined by:
%\begin{align}
%R(1,2)=\left(
%\begin{array}
%[c]{rrrr}%
%1 & -1 & 2 & -3\\
%0 & 0 & 1 & 0\\
%0 & 0 & 0 & 1
%\end{array}
%\right)  , &  R(2,1)=\left(
%\begin{array}
%[c]{rrr}%
%1 & 1 & 4\\
%-1 & 2 & 3\\
%0 & 1 & 0\\
%0 & 0 & 1
%\end{array}
%\right)  ,\nonumber\\
%R(2,3)=\left(
%\begin{array}
%[c]{rrr}%
%1 & -1 & 0\\
%0 & 1 & 1
%\end{array}
%\right)  , &  R(3,3)=\left(
%\begin{array}
%[c]{rr}%
%0 & 1\\
%2 & 10
%\end{array}
%\right)  ,\nonumber\\
%R(3,4)=\left(
%\begin{array}
%[c]{rr}%
%1 & 3
%\end{array}
%\right)  , &  R(4,1)=\left(
%\begin{array}
%[c]{r}%
%1\\
%-1\\
%0\\
%2
%\end{array}
%\right)  .\nonumber
%\end{align}
\end{itemize}
\begin{equation}
\begin{array}
[c]{ll}
R(1,2)=\left(
\begin{array}
[c]{rrrr}%
1 & -1 & 2 & -3\\
0 & 0 & 1 & 0\\
0 & 0 & 0 & 1
\end{array}
\right)  ,& R(2,1)=\left(
\begin{array}
[c]{rrr}%
1 & 1 & 0\\
-1 & 2 & 0\\
0 & 1 & 0\\
0 & 0 & 1
\end{array}
\right)  ,\nonumber\\
R(2,3)=\left(
\begin{array}
[c]{rrr}%
1 & -1 & 0\\
0 & 1 & 1
\end{array}
\right), &
R(2,5)=\left(
\begin{array}
[c]{rrr}%
1 & 0 & 0  \\
0 & 1 & 0 \\
0 & 0 & 1
\end{array}
\right),
\nonumber\\
 R(3,3)=\left(
\begin{array}
[c]{rr}%
1& 1\\
0 & 10
\end{array}
\right)  ,& %\nonumber\\%&
%\end{array}
%\end{equation}
%\begin{equation}
%\begin{array}
%[c]{ll}
R(3,6)=\left(
\begin{array}
[c]{rr}%
1 & 0\\
0 & 1
\end{array}
\right),\nonumber\\
R(4,1)=(
\begin{array}
[c]{rrrr}%
1& %\\
-1& %\\
0& %\\
2
\end{array}
)', &
R(4,2)=(
\begin{array}
[c]{rrr}%
1& %\\
1& %\\
1
\end{array}
)', \nonumber\\
R(5,4)=(
\begin{array}
[c]{rrr}%
1& %\\
1& %\\
1
\end{array}
),%\left(
%\begin{array}
%[c]{rrr}%
%1& %\\
%1& %\\
%1
%\end{array}
%\right) 
&
R(5,6)=
\left(
\begin{array}
[c]{rrr}%
2 & 1 & 0\\%\\
0 & 10 & 0
\end{array}
\right), 
\nonumber\\
R(6,5)=
\left(
\begin{array}
[c]{rr}%
1 & 0 \\
0 & 10 \\
0 & 10
\end{array}
\right). 
&.
\nonumber
% & C_{3}=\left(
%\begin{array}
%[c]{cc}%
%1 & 0
%\end{array}
%\right)  ,\nonumber\\
%A_{4}=3, &  B_{4}=1, & C_{4}=1;\nonumber
\end{array}
\end{equation}

\begin{figure}[h]
\begin{center}
\begin{tikzpicture}[->,shorten >=1pt,%
auto,node distance=2.75cm,semithick,
inner sep=1.5pt,bend angle=23]
%\draw[help lines] (3,3) grid (3,3);
\node[state] (q_1) {$1$};
\node[state] (q_2) [ right of=q_1] {$2$};
\node[state] (q_3) [ right of=q_2] {$3$};
\node[state] (q_4) [ below of=q_1] {$4$};
\node[state] (q_5) [ right of=q_4] {$5$};
\node[state] (q_6) [ right of=q_5] {$6$};
\path[->]
(q_1)
edge [bend left] node {} (q_2)
(q_2)
edge node {} (q_3)
edge [bend left] node {} (q_1)
edge node {} (q_5)
(q_3)
edge [loop right] node {} (q_3)
edge node {} (q_6)
(q_4)
edge node {} (q_1)
edge node {} (q_2)
(q_5)
edge node {} (q_4)
edge [bend left] node {} (q_6)
(q_6)
edge [bend left] node {} (q_5);
\end{tikzpicture}
\end{center}
\caption{Finite State Machine associated with the linear switching system
$\mathcal{H}$.}%
\label{fig1}%
\end{figure}
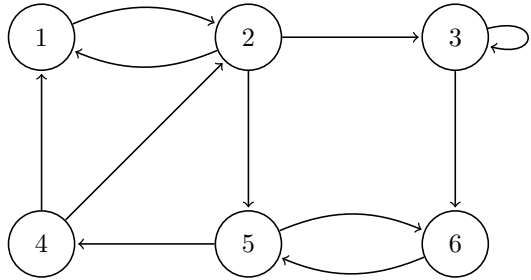

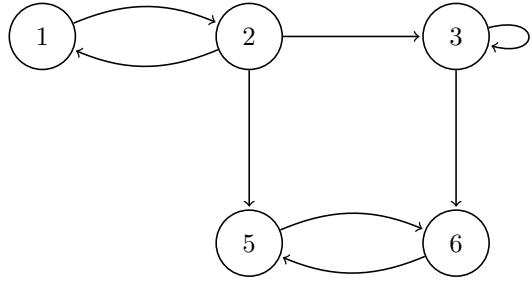
\begin{figure}[h]
\begin{center}
\begin{tikzpicture}[->,shorten >=1pt,%
auto,node distance=2.75cm,semithick,
inner sep=1.5pt,bend angle=23]
%\draw[help lines] (3,3) grid (3,3);
\node[state] (q_1) {$1$};
\node[state] (q_2) [ right of=q_1] {$2$};
\node[state] (q_3) [ right of=q_2] {$3$};
%\node[state] (q_4) [ below of=q_1] {$4$};
\node[state] (q_5) [ below of=q_2] {$5$};
\node[state] (q_6) [ right of=q_5] {$6$};
\path[->]
(q_1)
edge [bend left] node {} (q_2)
(q_2)
edge node {} (q_3)
edge [bend left] node {} (q_1)
edge node {} (q_5)
(q_3)
edge [loop right] node {} (q_3)
edge node {} (q_6)
%(q_4)
%edge node {} (q_1)
%edge node {} (q_2)
(q_5)
%edge node {} (q_4)
edge [bend left] node {} (q_6)
(q_6)
edge [bend left] node {} (q_5);
\end{tikzpicture}
\end{center}
\caption{Finite State Machine associated with the linear switching system
$\left.  \mathcal{H}^{\prime}\right\vert _{\widehat{Q}}$, with $\widehat
{Q}=\{1,2,3,5,6\}$.}%
\label{fig2}%
\end{figure}

The Finite State Machine associated with system $\mathcal{H}$ is depicted in
Figure \ref{fig1}.
Let us analyze observability and detectability properties of the linear switching system $\mathcal{H}$. The linear systems $S(i)$
associated with discrete states $i=1,2,3,5,6$ are detectable but not observable and therefore we conclude that $\mathcal{H}$ is
not observable. We now check detectability of $\mathcal{H}$. For this purpose, we apply Theorem \ref{thdetect}. We start by
checking condition i). The Markov parameters associated with systems $S(i)$, $i\in Q$ are given for any $k\in\mathbb{N}$ by:
%\[
\begin{equation}
\begin{array}
[c]{rrr}
%[c]{rrr}
C_{1}A_{1}^{k}B_{1}=1, & C_{2}A_{2}^{k}B_{2}=2^{k}, & C_{3}A_{3}^{k}B_{3}=0,\\
C_{4}A_{4}^{k}B_{4}=3^{k}, & C_{5}A_{5}^{k}B_{5}=4, & C_{6}A_{6}^{k}B_{6}=5^{k}.
\end{array}
\end{equation}
%\]
Hence, condition (\ref{condT}) is satisfied for $k=1$. Thus by Theorem \ref{Th_locobs}, the linear switching system $\mathcal{H}$
is location observable. We now check condition ii) of Theorem \ref{thdetect}. In this case $E^{\circlearrowright}=\{(3,3)\}$ and
\mbox{$Im(R(3,3)-I)\cap \ker(\mathcal{O}_{3})=\{0\}$}; thus condition ii) is satisfied. Finally, we check condition iii). Since
the linear system $S(4)$ is observable, by Proposition \ref{Th_dec_dis2} the switching system $\mathcal{H}^{\prime}$ associated
with $\mathcal{H}$ is detectable if and only if $\left.  \mathcal{H}^{\prime}\right\vert _{\widehat{Q}}$ with
$\widehat{Q}=\{1,2,3,5,6\}$, is detectable. The resulting linear switching system $\left. \mathcal{H}^{\prime}\right\vert _{\widehat
{Q}}$ is characterized by the Finite State Machine in Figure \ref{fig2}.\\
%\begin{figure}[hh]
%\begin{center}
%\begin{tikzpicture}[->,shorten >=1pt,%
%auto,node distance=2.75cm,semithick,
%inner sep=1.5pt,bend angle=23]
%%\draw[help lines] (3,3) grid (3,3);
%\node[state] (q_1) {$1$};
%\node[state] (q_2) [ right of=q_1] {$2$};
%\node[state] (q_3) [ right of=q_2] {$3$};
%\path[->]
%(q_1)
%edge [bend left] node {} (q_2)
%(q_2)
%edge node {} (q_3)
%edge [bend left] node {} (q_1)
%(q_3)
%edge [loop right] node {} (q_3);
%%edge node {} (q_4)
%%(q_4)
%%edge node {} (q_1)
%\end{tikzpicture}
%\end{center}
%\caption{Finite State Machine associated with the linear switching system
%$\left.  \mathcal{H}^{\prime}\right\vert _{\widehat{Q}}$, with $\widehat
%{Q}=\{1,2,3\}$.}%
%\label{fig2}%
%\end{figure}
We can now introduce the $GLSw$--system $\mathcal{H}_{0}$ of (\ref{S0}) associated with $\left.  \mathcal{H}^{\prime}\right\vert _{\widehat{Q}}$:
\begin{equation}
\mathcal{H}_{0}=\left(  \Xi_{0},S_{0},E,G_{0},R_{0}\right)  ,
\label{S0example}%
\end{equation}
where:

\begin{itemize}
\item $\Xi_{0}=\left(  \{1\}\times\mathbb{R}^{2} \right)  \cup\left(
\{2\}\times\mathbb{R}^{2} \right)  \cup\left(  \{3\}\times\mathbb{R}\right)  \cup\left(  \{5\}\times\mathbb{R}^{2}\right)  \cup\left(  \{6\}\times\mathbb{R}\right)  $;

\item $S_{0}(i)$ is described for any $i\in\widehat{Q}$ by dynamics
\mbox{$\dot {z}(t)=A_{i}^{(22)}z(t)$}, where:
\begin{align}
A_{1}^{(22)}=\left(
\begin{array}
[c]{rr}%
-2 & 1\nonumber\\
1 & -2
\end{array}
\right)  , & A_{2}^{(22)}=\left(
\begin{array}
[c]{rr}%
-1 & 1\\
1 & -2
\end{array}
\right)  & ,A_{3}^{(22)}=-1,\\
A_{5}^{(22)}=\left(
\begin{array}
[c]{rr}%
-1 & 0\\
0 & -2
\end{array}
\right)  , & A_{6}^{(22)}=-3\nonumber
\end{align}

\item $E=\{(1,2),(2,1),(2,3),(2,5),(3,3),(3,6),(5,6),(6,5)\}$;

\item $G_{0}(i,h)=\ker(R^{(12)}(i,h))$ for any $(i,h)\in E$, where:
\begin{equation}
\begin{array}
[c]{ll}
R^{(12)}(1,2)=(
\begin{array}
[c]{rr}%
2 & -3
\end{array}
),  &  R^{(12)}(2,1)=\left(
\begin{array}
[c]{rr}%
1 & 0\\
2 & 0
\end{array}
\right),  \\%&
R^{(12)}(2,3)=(
\begin{array}
[c]{rr}%
-1 & 0
\end{array}
),  &  R^{(12)}(2,5)=(
\begin{array}
[c]{rr}%
0 & 0
\end{array}
),\\
R^{(12)}(3,3)=1, & R^{(12)}(3,6)=0,\\
R^{(12)}(5,6)=(
\begin{array}
[c]{rr}
1 & 0
\end{array}
), & R^{(12)}(6,5)=0;
\end{array}
\nonumber
\end{equation}

\item $R_{0}(i,h)=R^{(22)}(i,h)$ for any $(i,h)\in E$, where:
\begin{equation}
\begin{array}
[c]{ll}
R^{(22)}(1,2)=\left(
\begin{array}
[c]{rr}%
1 & 0\\
0 & 1
\end{array}
\right),  &  R^{(22)}(2,1)=\left(
\begin{array}
[c]{rr}%
1 & 0\\
0 & 1
\end{array}
\right), \\%&
R^{(22)}(2,3)=(
\begin{array}
[c]{rr}%
1 & 1
\end{array}
),  &  
R^{(22)}(2,5)=\left(
\begin{array}
[c]{rr}%
1 & 0\\
0 & 1
\end{array}
\right),\\
R^{(22)}(3,3)=10, &
R^{(22)}(3,6)=1, \\
R^{(22)}(5,6)=(
\begin{array}
[c]{rr}%
10 & 0
\end{array}
), &
R^{(22)}(6,5)=\left(
\begin{array}
[c]{rr}%
10 \\ 10
\end{array}
\right).
\end{array}
\nonumber
\end{equation}
\end{itemize}

The Finite State Machine associated with $\mathcal{H}_{0}$ (Figure \ref{fig2}) is composed by three strongly connected components\footnote{We recall that a strongly connected component of a FSM is a FSM, with a path between any two discrete
states.},
i.e. one involving the discrete states $1,2\in\hat{Q}$, the one involving the discrete state $3\in\hat{Q}$ and the other involving the discrete states $5,6\in\hat{Q}$. It is well--known
that the asymptotic stability of a switching system can be assessed by studying this property in each strongly connected
component. More precisely, $\mathcal{H}_{0}$ is
asymptotically stable if and only if switching system $\left.  \mathcal{H}%
_{0}\right\vert _{Q_{1}}$, with \mbox{$Q_{1}=\{1,2\}$}, switching system $\left.  \mathcal{H}%
_{0}\right\vert _{Q_{2}}$, with \mbox{$Q_{2}=\{3\}$}, and switching system $\left.
\mathcal{H}_{0}\right\vert _{Q_{3}}$, with \mbox{$Q_{3}=\{5,6\}$} are asymptotically
stable.\newline 
We first consider $\left.  \mathcal{H}%
_{0}\right\vert _{Q_{1}}$. %By Proposition \ref{Th_detect} the $GLSw$--system $\left.  \mathcal{H}%
%_{0}\right\vert _{Q_{1}}$ is asymptotically stable if the corresponding abstraction $\mathcal{H}_{1}$, defined in (\ref{abs}) is asymptotically stable. 
We recall from \cite{LiberzonBook} that an
autonomous $GLSw$--system (with identity reset map) is asymptotically
stable if it admits a common Lyapunov function $V$. By defining for
any $x\in\mathbb{R}^{2}$ the function $V(x)=x^{\prime}Px$ with $P=I$, we obtain:
\begin{align}
(A_{1}^{(22)})^{\prime}\,P+P\,A_{1}^{(22)}\leq-\mathcal{Q}, \hspace{0.2cm}&%\nonumber\\
(A_{2}^{(22)})^{\prime}\,P+P\,A_{2}^{(22)}\leq-\mathcal{Q},\nonumber
\end{align}
where:
\[
\mathcal{Q}=\left[
\begin{array}
[c]{rr}%
-2 & 2\\
2 & -4
\end{array}
\right]  \geq0.
\]
Hence $V$ is a common Lyapunov function for sub--systems $S_{0}(1)$ and $S_{0}(2)$ of $\left.
\mathcal{H}_{0}\right\vert _{Q_{1}}$ and by Theorem 2.1 in \cite{LiberzonBook} we conclude that $\left.  \mathcal{H}_{0}\right\vert _{Q_{1}}$ is
asymptotically stable\footnote{Dynamical matrices $A_{1}^{(22)}$ and $A_{2}%
^{(22)}$ have been taken from \cite{Narendra.tac1994}.}.
\newline Let us consider $\left. \mathcal{H}_{0}\right\vert
_{Q_{2}}$. The $GLSw$--system $\left.
\mathcal{H}_{0}\right\vert _{Q_{2}}$ is characterized by dynamical matrix $A_{3}^{(22)}=-1$, reset matrix $R_{2}(3,3)=0$ and guard $G_{0}(3,3)=\{0\}$ and hence it is asymptotically stable.\\
Let us now consider $\left. \mathcal{H}_{0}\right\vert
_{Q_{3}}$ and let us apply Proposition \ref{Th_detect} to investigate stability properties of $\left. \mathcal{H}_{0}\right\vert
_{Q_{3}}$. It is readily seen that the abstraction $\mathcal{H}_{1}$ of (\ref{abs}), that corresponds to $\left. \mathcal{H}_{0}\right\vert
_{Q_{3}}$ is unstable. Let us now consider the abstraction $\mathcal{H}_{2}$ of $\left. \mathcal{H}_{0}\right\vert _{Q_{3}}$. The reset map $R_{2}(e)$ with $e=(5,6)$ associated to $\mathcal{H}_{2}$ is given by:
\[
R_{2}(5,6) = R^{(22)}(5,6)\pi_{\ker(R^{(12)}(5,6))}=(
\begin{array}
[c]{rr}
0 & 0
\end{array}
).
\]
%\end{equation}
Therefore, since dynamical matrices $A^{(22)}_{5}$ and $A^{(22)}_{6}$ are Hurwitz it is easy to see that the $LSw$--system  $\mathcal{H}_{2}$ is asymptotically stable. Thus, by Proposition \ref{Th_detect} also $\left. \mathcal{H}_{0}\right\vert _{Q_{3}}$ is asymptotically stable.\\
We conclude that the switching system $\mathcal{H}_{0}$ is asymptotically stable and therefore
condition iii) of Theorem \ref{thdetect} is satisfied. Hence, by Theorem \ref{thdetect}, the linear switching system
$\mathcal{H}$ is detectable.

\section{Conclusions}

We addressed observability and detectability of linear switching systems.
%We
%proposed a definition of observability, and a weaker notion of detectability,
%related to the possibility of reconstructing the system state.
We derived a computable necessary and sufficient condition for a switching system to be observable. Further, we derived a Kalman
decomposition of the switching system, which reduces detectability of linear switching systems to asymptotic stability of
suitable linear switching systems with guards associated with the original systems. The study of detectability is a fundamental
step towards the design of a hybrid observer. In fact, by Definition \ref{def_detect}, a necessary condition for the existence of
a hybrid observer for a $LSw$--system $\mathcal{H}$ is that $\mathcal{H}$ is detectable. On the other hand, as shown in Section
3, observability of $\mathcal{H}$ implies the existence of an algorithm that reconstructs the current hybrid state; in
particular, the combination of (\ref{DiscreteObs}) and (\ref{ContinuousObs}) can be thought of as a hybrid observer. However,
such an observer requires an infinite precision in the computation of the vector $Y^{(n)}(t)$. Further work will identify appropriate conditions on linear switching systems, for the existence and design of hybrid observers.

\bibliographystyle{plain}
\bibliography{byblioAutom}

\section{Appendix: Proof of Lemma \ref{Th_nonempty}}

We first need two preliminary technical lemmas.

\begin{lemma}
\label{L1} If condition (\ref{condT}) is satisfied then for any $(i,h)\in
\hat{J}$,
%
%Given $\left(  i,h\right)  \in\hat{J}$, if there exists
%$k\in\mathbb{N}$, $k<n_{i}+n_{h}$, such that $C_{i}A_{i}^{k}B_{i}\neq
%C_{h}A_{h}^{k}B_{h}$, then
$B_{ih}^{-1}\left(  \mathcal{V}_{ih}\right)  \neq\mathbb{R}^{m}$.
\end{lemma}

\begin{proof}
By contradiction, suppose that $B_{ih}^{-1}\left(  \mathcal{V}_{ih}\right)
=\mathbb{R}^{m}$ for some $(i,h)\in\hat{J}$. Then $Im\left(  B_{ih}\right)
\subseteq\mathcal{V}_{ih} $ and by (\ref{inv}), $A_{ih}\mathcal{V}_{ih}%
\subseteq\mathcal{ V}_{ih}+Im\left(  B_{ih}\right)  \subseteq\mathcal{ V}_{ih} $,
%and therefore the above condition implies that $A_{ih}\mathcal{V}%
%_{ih}\mathcal{\subset V}_{ih}$,
i.e. $\mathcal{V}_{ih}$ is $A_{ih}-$invariant and contains $Im\left(
B_{ih}\right)  $. Since the minimal $A_{ih}-$invariant subspace containing
$Im\left(  B_{ih}\right)  $ is $Im(
\begin{array}
[c]{cccc}%
B_{ih} & A_{ih}B_{ih} & \ldots & A_{ih}^{n-1}B_{ih}%
\end{array}
)  $, with $n=n_{i}+n_{h}$, then
%. Therefore%
%\[
$
Im(
\begin{array}
[c]{cccc}%
B_{ih} & A_{ih}B_{ih} & \ldots & A_{ih}^{n-1}B_{ih}%
\end{array}
)  \subseteq\mathcal{V}_{ih}\subseteq\ker\left(  C_{ih}\right)  ,
$
%\]
which implies $
%\[
C_{ih}(
\begin{array}
[c]{cccc}%
B_{ih} & A_{ih}B_{ih} & \ldots & A_{ih}^{n-1}B_{ih}%
\end{array}
)  =0
%\]
$. Thus condition (\ref{condT}) is not satisfied and hence a contradiction holds.
\end{proof}

\begin{lemma}
\label{Lemma_lambda}Let $\left\{  M_{i}\in\mathbb{R}^{m\times nT},i\in
Q\right\}  $ be a family of nonzero matrices. There exists $z\in\mathbb{R}^{n}
$ and $\lambda\in\mathbb{R}$ such that $M_{i}\mathbf{z}\neq0$, $\forall i\in
Q$, where
%\begin{equation}
%\mathbf{z}\mathbf{=}\left(
%\begin{array}
%[c]{c}%
%z\\
%\lambda z\\
%\lambda^{2}z\\
%\vdots\\
%\lambda^{T-1}z
%\end{array}
%\right)  .%\label{zetabold}%
%\end{equation}%
\begin{equation}
\mathbf{z^{\prime}}\mathbf{=}(
\begin{array}
[r]{rrrrr}%
z^{\prime} & \lambda z^{\prime} & \lambda^{2}z^{\prime} & \hdots &
\lambda^{T-1}z^{\prime}%
\end{array}
)  ^{\prime}. \label{zetabold}%
\end{equation}

\end{lemma}

\begin{proof}
By setting
%\[
%\begin{array}
%[c]{l}%
%M_{i}=(
%\begin{array}
%[c]{ccccc}%
%M_{i0} & M_{i1} & M_{i2} & \ldots & M_{iT-1}%
%\end{array}
%)  ,\\
%M_{i}(\lambda)=M_{i0}+\lambda M_{i1}+\lambda^{2}M_{i2}+\ldots+\lambda
%^{T-1}M_{iT-1},
%\end{array}
%\]
$
M_{i}=(
\begin{array}
[c]{cccc}%
M_{i0} & M_{i1} & \ldots & M_{iT-1}%
\end{array}
)$ and $M_{i}(\lambda)=M_{i0}+\lambda M_{i1}+\lambda^{2}M_{i2}+\ldots+\lambda
^{T-1}M_{iT-1}$,
with $M_{ij}\in\mathbb{R}^{m\times n}$, for any $z\in\mathbb{R}^{n}$,
$M_{i}\mathbf{z}=M_{i}(\lambda)z$. Given $i\in Q$, since $M_{i}\neq0$, there
are a finite number of values $\theta$ such that $M_{i}(\theta)=0$. Choose
$\lambda$ such that $M_{i}(\lambda)\neq0$, $\forall i\in Q$. Then there exists
$z\notin\bigcup_{i\in Q}ker\left(  M_{i}(\lambda)\right)  $ which implies
$M_{i}\mathbf{z}\neq0$, $\forall i\in Q$.
\end{proof}

\bigskip

We now give the proof of Lemma \ref{Th_nonempty}.

\begin{proof}
By contradiction, suppose that the set $\mathcal{U}^{\ast}$ is empty and let be
$n=n_{i}+n_{h}$. Then
\begin{equation}
\begin{array}
[c]{l}%
\forall u\in\mathcal{U}\text{, }\exists t^{\prime},t^{\prime\prime}%
\in\mathbb{R},\exists\left(  i,h\right)
\in\hat{J}\text{ and }\widetilde{u}\in\mathcal{U}_{ih}\text{ }s.t.\\
%\text{}
u(t)=\widetilde{u}(t)\text{,\ }\forall t\in\left[  t^{\prime},t^{\prime
\prime}\right]  . %
\end{array}
\label{uno}
\end{equation}
Let $V^{ih}$ be the set of smooth functions $v:\mathbb{R\rightarrow}%
B_{ih}^{-1}\left(  \mathcal{V}_{ih}\right)  $ and let $\widehat{\mathcal{U}%
}\subset\mathcal{U}$ be the set of smooth, not identically zero functions. By
definition of $\mathcal{U}_{ih}$, condition $\left(  \ref{uno}\right)  $
implies:%
\begin{equation}%
\begin{array}
[c]{l}%
\forall u\in\widehat{\mathcal{U}}\text{, }\exists t^{\prime},t^{\prime\prime
}\in\mathbb{R},\exists\left(  i,h\right)
\in\hat{J}\text{, }\overline{z}\in\mathcal{V}_{ih}%\text{ }%\\
\text{ and }v_{ih}\in V^{ih}\\
s.t. \text{ }u(t)=K_{ih}z(t)+v_{ih}(t)\text{,\ }\forall t\in\left[  t^{\prime}%
,t^{\prime\prime}\right],
\end{array}
\label{due}%
\end{equation}
where $\dot{z}(t)=\widehat{A}_{ih}z(t)+B_{ih}v_{ih}(t),$ $\widehat{A}%
_{ih}=A_{ih}+B_{ih}K_{ih}$ and $z(t^{\prime})=\overline{z}\in\mathcal{V}_{ih}
$. Condition $\left(  \ref{due}\right)  $ implies:%
\begin{align}
\forall u  &  \in\widehat{\mathcal{U}}\text{, }\exists t^{\prime}\in
\mathbb{R}\text{, }\exists\left(  i,h\right)  \in\hat{J}\text{ s.t. }%
\forall\bar{N}\geq0\label{cond}\\
\left(
\begin{array}
[c]{c}%
u(t^{\prime})\\
\dot{u}(t^{\prime})\\
%\ddot{u}(t^{\prime})\\
\ldots\\
u^{(\bar{N})}(t^{\prime})
\end{array}
\right)   &  \in\mathbf{M}_{ih}^{\bar{N}}\mathcal{V}_{ih}+\mathbf{F}%
_{ih}^{\bar{N}}\left(  F_{ih}\times F_{ih}\times\ldots\times F_{ih}\right)
,\nonumber
\end{align}
where $F_{ih}=B_{ih}^{-1}\left(  \mathcal{V}_{ih}\right)$ and
\begin{eqnarray}
%\begin{align*}
%\begin{array}
%[c]{cc}%
& \mathbf{M}%
_{ih}^{\bar{N}}=\left(
\begin{array}
[c]{c}%
K_{ih}\\
K_{ih}\widehat{A}_{ih}\\
%K_{ih}\widehat{A}_{ih}^{2}\\
\ldots\\
K_{ih}\widehat{A}_{ih}^{\bar{N}}%
\end{array}
\right)  \in\mathbb{R}^{m\left(  \bar{N}+1\right)  \times\bar{N}},\nonumber\\ %&
%\end{array}
%&
%\hspace{0.5cm}
&
\mathbf{F}_{ih}^{\bar{N}} =\left(
\begin{array}
[c]{cccc}%
I & 0 & \ldots & 0\\
K_{ih}B_{ih} & I & \ldots & 0\\
\ldots & \ldots & \ldots & 0\\
K_{ih}\widehat{A}_{ih}^{\bar{N}-1}B_{ih} & K_{ih}\widehat{A}_{ih}^{\bar{N}%
-2}B_{ih} & \ldots & I
\end{array}
\right)  \in\mathbb{R}^{m\bar{N}\times m\bar{N}}.\nonumber
%\end{align*}
\end{eqnarray}

The matrix $\mathbf{F}_{ih}^{\bar{N}}$ is nonsingular. By setting $\dim
(F_{ih})=\nu$, one obtains:
\[
\dim(\mathbf{F}_{ih}^{\bar{N}}\left(  F_{ih}\times F_{ih}\times\ldots\times
F_{ih}\right)  )=\nu\left(  \bar{N}+1\right)  ,
\]
and since (\ref{condT}) holds,\ $\dim(\mathbf{M}_{ih}^{\bar{N}}\mathcal{V}%
_{ih})<n$; thus
\[
\dim(\mathbf{M}_{ih}^{\bar{N}}\mathcal{V}_{ih}+\mathbf{F}_{ih}^{\bar{N}%
}\left(  F_{ih}\times F_{ih}\times\ldots\times F_{ih}\right)  )\leq\nu\left(
\bar{N}+1\right)  +n.
\]
Therefore since by Lemma \ref{L1}, $\nu<m$, we obtain that $\nu\left(  \bar{N}+1\right)  +n<m\left(  \bar{N}+1\right)  $ for any $\bar{N}>\frac{n}{m-\nu}-1$; thus
$\mathbf{M}_{ih}^{\bar{N}}\mathcal{V}_{ih}+\mathbf{F}_{ih}^{\bar{N}%
}\left(  F_{ih}\times F_{ih}\times\ldots\times F_{ih}\right)  $ is a proper
subspace of $\mathbb{R}^{m\left(  \bar{N}+1\right)  }$. Hence there exists a
sufficiently large $\bar{N}$ such that the set $\mathbf{M}_{ih}^{\bar{N}%
}\mathcal{V}_{ih}+\mathbf{F}_{ih}^{\bar{N}}\left(  F_{ih}\times F_{ih}%
\times\ldots\times F_{ih}\right)  $ is a proper subspace of $\mathbb{R}%
^{m\left(  \bar{N}+1\right)  }$ for any $(i,h)\in\hat{J}$. Given some
$z\in\mathbb{R}^{m}$ and $\lambda\in\mathbb{R}$ let be $\mathbf{u}%
(t)=z\exp\left(  \lambda t\right)  \in\widehat{\mathcal{U}}$. It follows that:%
\[
\left(
\begin{array}
[c]{c}%
\mathbf{u}(t)\\
\mathbf{\dot{u}}(t)\\
%\mathbf{\ddot{u}}(t)\\
\ldots\\
\mathbf{u}^{(\bar{N})}(t)
\end{array}
\right)  =\left(
\begin{array}
[c]{c}%
z\\
\lambda z\\
%\lambda^{2}z\\
\vdots\\
\lambda^{T-1}z
\end{array}
\right)  \exp\left(  \lambda t\right)  .
\]
Set $\mathbf{M}_{ih}^{\bar{N}}\mathcal{V}_{ih}+\mathbf{F}_{ih}^{\bar{N}%
}\left(  F_{ih}\times F_{ih}\times\ldots\times F_{ih}\right)  =ker \left(
G_{ih}\right)  $, for some matrix $G_{ih}$. By Lemma
\ref{Lemma_lambda} there exist $z$ and $\lambda$ such that $G_{ih}%
\mathbf{z}\neq0$, $\forall\left(  i,h\right)  \in\hat{J}$ where $\mathbf{z}$
is as in (\ref{zetabold}). This implies that the vector
\[
\left(
\begin{array}
[c]{c}%
\mathbf{u}(t)\\
\mathbf{\dot{u}}(t)\\
%\mathbf{\ddot{u}}(t)\\
\ldots\\
\mathbf{u}^{(\bar{N})}(t)
\end{array}
\right)  =\mathbf{z}\exp(\lambda t)
\]
does not belong to $\mathbf{M}_{ih}^{\bar{N}}%
\mathcal{V}_{ih}+\mathbf{F}_{ih}^{\bar{N}}\left(  F_{ih}\times F_{ih}%
\times\ldots\times F_{ih}\right)$, for all $(i,h)\in
\hat{J}$ and $t\in\mathbb{R},$
and hence condition $\left(  \ref{cond}\right)  $ is false; thus the result follows.
\end{proof}

\end{document}